\documentclass[12pt]{amsart}
 \usepackage{tikz}
\usepackage{amssymb,amsmath}
\newtheorem{theorem}{Theorem} 

\numberwithin{equation}{section}

\def\di{\partial}
\def\dib{\bar\partial}
\numberwithin{equation}{section}
\def\simleq{\underset\sim<}
\def\simgeq{\underset\sim>}

\def\T{\text}
\def\1#1{\overline{#1}}
\def\2#1{\widetilde{#1}}
\def\3#1{\widehat{#1}}
\def\4#1{\mathbb{#1}}
\def\5#1{\frak{#1}}
\def\6#1{{\mathcal{#1}}}
\def\C{{\4C}}

\def\sumK{\underset{|K|=k-1}{{\sum}'}}

\def\sumij{\underset {ij=1,\dots,N}{{\sum}}}

\def\contrazione{\raisebox{1pt}{\,{\mbox{\tiny{$|\!\raisebox{-0.7pt}{\underline{\hphantom{X}}}$}}}\,}}

\def\phi{\varphi}

\newtheorem{Thm}{Theorem}
\newtheorem{Cor}[Thm]{Corollary}
\newtheorem{Pro}[Thm]{Proposition}
\newtheorem{Lem}[Thm]{Lemma}
\theoremstyle{definition}\newtheorem{Def}[Thm]{Definition}
\theoremstyle{remark}
\newtheorem{Rem}[Thm]{Remark}
\newtheorem{Exa}[Thm]{Example}

\def\Label#1{\label{#1}}
\def\bl{\begin{Lem}}
\def\el{\end{Lem}}
\def\bp{\begin{Pro}}
\def\ep{\end{Pro}}
\def\bt{\begin{Thm}}
\def\et{\end{Thm}}
\def\bc{\begin{Cor}}
\def\ec{\end{Cor}}
\def\bd{\begin{Def}}
\def\ed{\end{Def}}
\def\br{\begin{Rem}}
\def\er{\end{Rem}}
\def\be{\begin{Exa}}
\def\ee{\end{Exa}}
\def\bpf{\begin{proof}}
\def\epf{\end{proof}}
\def\ben{\begin{enumerate}}
\def\een{\end{enumerate}}

\def\1alpha{[\frac1\alpha]}
\def\T{\text}

\def\C{{\Bbb C}}

\numberwithin{equation}{section}
\def\T{\text}

\newcommand{\no}[1]{\|{#1}\|}
\newcommand{\NO}[1]{\|{#1}\|^2}

\newtheorem{definition}[theorem]{Definition }
\newtheorem{lemma}[theorem]{Lemma  }
\newtheorem{proposition}[theorem]{Proposition  }
\newtheorem{corollary}[theorem]{Corollary }
\newtheorem{example}[theorem]{\it Example }

\begin{document}
\title[Extension of $L^2$ holomorphic functions]
{Extension of $L^2$ holomorphic functions}
\author[L.~Baracco ]{Luca Baracco}
\address{Dipartimento di Matematica, Universit\`a di Padova, via 
Trieste 63, 35121 Padova, Italy}
\email{baracco@math.unipd.it}  
%\author[ T.V.~Khanh and G.~Zampieri ]
%{Tran Vu Khanh and Giuseppe Zampieri}
%\address{Dipartimento di Matematica, Universit\`a di Padova, via 
%Trieste 63, 35121 Padova, Italy}
%\email{khanh@math.unipd.it,
%zampieri@math.unipd.it}
\maketitle
\begin{abstract}
The purpose of this note is to show that the $\dib$-estimate which is needed in the Ohsawa-Takegoshi Extension Theorem \cite{O88} is a direct consequence of the H\"ormander-Kohn-Morrey weigthed inequality. In this inequality, the Donnelly-Fefferman argument is not required and a single 1-parameter family of non-singular weights is used.
%%%
%and does not require the Donnelly-Fefferman argument nor the use of singular weights. 
%in which the gradient of the weight $\di\phi$ needs not to be ``absorbed"
% by the Levi form $\di\dib\phi$. (This property is only implicitly used when assuming the positivity of $\di\dib\phi$.) 
This paper is the furtherst step of a great deal of work devoted to the  simplification of the original proof of Ohsawa-Takegoshi Theorem; among other papers on the subject, we mention \cite{B96} and \cite{S96} which are based on ``twisted" basic estimates and, in recent time, \cite{C11} and \cite{S11}.

\end{abstract}

Through an approximation argument, the extension theorem relies on Theorem~\ref{t1.1} below.
\bt
\Label{t1.1}
Let $D\subset\subset\C^n$ be a bounded smooth pseudoconvex domain with diameter $\le1$, $\psi$ a plurisubharmonic function on $D$, $\alpha$ a $\dib$-closed form in $L^2_\psi(D)$ of degree $\ge 1$ such that $\alpha_J=0$ for $1\notin J$ and supp$\,\alpha\subset\{z\,:\,|z_1|<\delta\}$. Then there is a solution $u=u_\delta\in L^2_\psi$ to the problem
\begin{equation}
\Label{1.1}
\begin{cases}
\dib u=\alpha,
\\
\no{z_1 u}_\psi\le c\delta\no{\alpha}_\psi,\quad\T{for $c$ independent of $\delta$, $\psi$ and $D$}.
\end{cases}
\end{equation}
\et
%\br
% more involved argument, which uses the selfboundedness of the gradient of the weights, gives the more accurate estimate
%begin{equation*}
%\no{(-\log|z_1|)^{-\frac12-\frac\eta2}u}_\psi\le c\no{\alpha}_\psi,
%\end{equation*}
%for any $\eta>0$ and for $c=c_\eta$ uniform in $\delta$. However, \eqref{1.1} suffices for extension.
%\er
\bpf
We first assume $\psi=0$.
We set
\begin{equation}
\Label{1.2}
\phi(=\phi^\delta)=-\log(-\log(|z_1|^2+\delta^2)),\qquad a(=a^\delta)=-\log(|z_1|^2+\delta^2);
\end{equation}
they are related by $a=e^{-\phi}$.  Their  key properties are
\begin{equation}
\Label{1.a}
\phi_1^\delta=\frac {-\bar z_1}{\log(|z_1|^2+\delta^2)(|z_1|^2+\delta^2)},\quad \phi_{1\bar 1}^\delta=\frac{|z_1|^2-\delta^2\log(|z_1|^2+\delta^2)}{\log^2(|z_1|^2+\delta^2)(|z_1|^2+\delta^2))^2}.
\end{equation}
This readily implies
\begin{equation}
\Label{1.b}
-a_{1\bar1}\ge0,\quad-a_{1\bar1}\simgeq \frac1{\delta^2}\T{ for $|z_1|<\delta$}.
\end{equation}
The  contraction of the gradient with a multivector $v$ is defined by $\di a\contrazione v=(a_i)\contrazione v=\sum_ia_iv_{iK}$ and the action of the Levi form is $\di\dib a(v,v)=(a_{i\bar j})(v,v)=\sumK\sumij a_{i\bar j}v_{iK}\bar v_{jK}$. 
From the identity  $\dib^*_\phi=\dib^*+\di\phi\contrazione$, we get the estimate
\begin{equation}
\Label{1.4}
\NO{\dib^*_\phi v}_\phi\le \NO{\dib^*v}_\phi+\NO{\di\phi\contrazione v}_\phi+2\Big|\int_De^{-\phi}(\di\phi\contrazione v)\cdot (\dib^*v)dV\Big|.
\end{equation}
We also have
\begin{equation}
\Label{1.5}
\begin{split}
2\Big|\int_De^{-\phi}\di\phi\contrazione v\,\cdot\,\dib^*v\,dV\Big|&\underset{\T{\eqref{1.2}}}=2\Big|\int_D\di a\contrazione v \cdot \dib^* v dV\Big|
\\
&\underset{\T{Cauchy-Schwarz}}\le \no{v}_0^2+\NO{|\di a|\dib^* v}_0
\\
&\underset{\T{Pseudoconvexity}}\simleq c(\NO{\dib v}_0+\NO{\dib^*v}_0+\NO{|\di a|\dib^* v}_0),
\end{split}
\end{equation}
where $c$ only depends on the diameter of $D$. Recall that $a=e^{-\phi}$ and $\no{\cdot}_\phi=\no{a^{\frac12}\cdot}_0$. We may then conclude
\begin{equation}
\Label{1.1bis}
\begin{split}
\int_D-(a_{i\bar j})(v,v) dV&=\int e^{-\phi}\phi_{i\bar j}(v,v) dV-\NO{\di\phi\contrazione u}_\phi
\\
&\underset{\T{basic}}\le \NO{\dib v}_\phi+\NO{\dib^*_\phi v}_\phi-\NO{\di\phi\contrazione u}_\phi
\\
&\underset{\T{\eqref{1.4}}}\le \NO{\dib v}_\phi+\NO{\dib^* v}_\phi +
2\Big|\int_De^{-\phi}\di\phi\contrazione v\,\cdot\,\dib^*v\,dV\Big|
\\
&\le c(\NO{(1+a^{\frac12})\dib v}_0+\NO{(1+a^{\frac12}+|\di a|)\dib^*v}_0).
\end{split}
\end{equation}
With \eqref{1.1bis} in our hands, 
%we define a functional in $L^2$. For $v\in \cap D_{\dib^*}$, we set
%\begin{equation}
%\Label{2.10}
%(1+a^{\frac12}+|\di a|)\dib^*v\mapsto
%\begin{cases}
%(v,\alpha)_0&\T{if $v\in \ker\dib$}
%\\
%0&\T{if $v\in(\ker\dib)^\perp$}.
%\end{cases}
%\end{equation}
%We have for $v\in\ker\dib$
we  define a functional in $L^2$ by putting
\begin{equation}
\Label{2.10}
(1+a^{\frac12}+|\di a|)\dib^*v\mapsto
(v,\alpha)_0\quad\T{for $v\in D_{\dib^*}$}. 
\end{equation}
If $v\in (\ker\dib)^\perp$, we have $(v,\alpha)_0=0$. If, instead, $v\in \ker\dib$, then 
\begin{equation}
\Label{2.11}
\begin{split}
|(v,\alpha)_0|&=\Big|\int_D(-a_{i\bar j})(v,(-a_{1\bar 1})^{-1}\alpha)dV\Big|
\\
&
\underset{\T{Cauchy-Schwarz}}\le \left(\int_D(-a_{i\bar j})(v,v)dV\right)^{\frac12}\no{(-a_{1\bar1})^{-\frac12}\alpha}_0
\\
&\underset{\T{\eqref{1.b}}}\le \Big(\int_D(-a_{i\bar j})(v,v)dV\Big)^{\frac12}(\delta\no{\alpha}_0)
\\
&\underset{\T{\eqref{1.1bis}}}\simleq \no{(1+a^{\frac12}+|\di a|)\dib^* v} \,\delta\no{\alpha}_0,
\end{split}
\end{equation}
where, in order to apply Cauchy-Schwarz, we have used that $(-a_{i\bar j})\ge0$.
Thus, by Riesz Representation Theorem, there exists $w\in L^2$ of norm $\no{w}_0\simleq\delta\no{\alpha}_0$ which represents \eqref{2.10}. 
Setting $u:=(1+a^{\frac12}+|\di a|)w$, we get $\dib u=\alpha$; also, from the trivial inequality $(1+a^{\frac12}+|\di a|)\le |z_1|^{-1}$, we conclude $\no{z_1u}_0\simleq\delta\no{\alpha}_0$.  (Note that we replace $(1+a^{\frac12}+|\di a|)$ by $|z_1|^{-1}$ in the conclusion but not in the proof, especially in \eqref{2.11}, because $|z_1|^{-1}$ is not in $L^2$.)
\vskip0.3cm
\noindent
When $\psi$ is smooth in $\bar D$, the proof above can be repeated verbatim, with the weight $\phi $ replaced by $\phi+\psi$ in the basic estimate \eqref{1.1bis} and with the Levi form $(\psi_{ij})$ dropped down because of its positivity. For a general plurisubharmonic $\psi$, we take a smooth  approximation of $\psi$ from above and an exhaustion of $D$ from inside. A weak limit of solutions on these subdomains yields the solution on $D$. For this approximation we need that $L^2_\psi\subset L^2_{\T{loc}}$. But this follows from the local boundedness of $\psi$ from above which is in turn a consequence of its upper semicontinuity.

\epf

\end{document}